\documentclass{article}
\usepackage{amssymb}
\usepackage{natbib}

\newtheorem{lemma}{Lemma}[section]
\newtheorem{theorem}[lemma]{Theorem}
\newtheorem{proposition}[lemma]{Proposition}
\newenvironment*{proof}{\noindent\textit{Proof}:~}{\vskip 2mm}

\def\field{F}                           % generic notation for a field        %
\def\sml#1{M(#1)}                       % the Paige loop constructed over #1  %
\def\aut#1{\mathrm{Aut}(#1)}            % the group of automorphisms of #1    %
\def\mlt#1{\mathrm{Mlt}(#1)}            % multiplication group of #1          %
\def\neutral{1}                         % generic notation for a neutral elm. %
\def\octo#1{O(#1)}                      % the split octonion algebra over #1  %
\def\to{\longrightarrow}                % arrow for a map                     %
\def\uocto#1{O(#1)^*}                   % unit elements in \octo{#1}          %
\def\net{\mathcal N}                    % generic 3-net                       %
\def\points{\mathcal P}                 % points of a generic 3-net           %
\def\lines{\mathcal L}                  % lines of a generic 3-net            %
\def\centre#1{Z(#1)}                    % center of #1                        %
\def\origin{1}                          % origin of a 3-net                   %
\def\bol#1{\sigma_{#1}}                 % bol reflection with axis #1         %
\def\comm#1#2{[#1,\,#2]}                % commutator of #1 and #2             %
\def\spn#1{\langle #1\rangle}           % subalgebra spanned by #1            %
\def\coll#1{\mathrm{Coll}(#1)}          % collineation group of a net #1      %
\def\sqmap#1{\varphi_{#1}}              % map of \net induced by an           %
\def\indmap#1{\widehat{#1}}             % the automorphism of \Gamma induced  %
\def\centralizer#1#2{C_{#1}(#2)}        % centralizer of #2 in #1             %
\def\chr#1{\mathrm{char}\ #1}           % characteristic of a field #1        %
\def\inn#1{\mathrm{Inn}(#1)}            % inner automorphisms of a group #1   %

\begin{document}
\title{
\Large{
\textbf{Automorphism Groups of Simple Moufang Loops over Perfect Fields}
}}

\author{\textsc{By G\'ABOR P.~NAGY}
\thanks{Supported by the ``J\'anos Bolyai Fellowship'' of the Hungarian
Academy of Sciences, and by the grants OTKA T029849 and FKFP
0063/2001.}
\\
\textit{SZTE Bolyai Institute}
\\
\textit{Aradi v\'ertan\'uk tere $1$, H-$6720$ Szeged, Hungary}
\\
\textit{e-mail}: \texttt{nagyg@math.u-szeged.hu}
\and
\textsc{PETR VOJT\v{E}CHOVSK\'Y}
\thanks{ Partially supported by the Grant Agency of Charles University, grant no.
$269$/$2001$/B-MAT/MFF, and by research assistantship at Iowa State University.}
\\
\textit{Department of Mathematics, Iowa State University}
\\
\textit{$400$ Carver Hall, Ames, Iowa, $50011$, USA}
\\
\textit{e-mail}: \texttt{petr@iastate.edu}
}

\maketitle

\begin{abstract}
Let $\field$ be a perfect field and $\sml{F}$ the nonassociative simple
Moufang loop consisting of the units in the (unique) split octonion algebra
$\octo{\field}$ modulo the center. Then $\aut{\sml{F}}$ is equal to
$G_2(\field) \rtimes \aut{\field}$. In particular, every automorphism of
$\sml{\field}$ is induced by a semilinear automorphism of $\octo{\field}$. The
proof combines results and methods from geometrical loop theory, groups of Lie
type and composition algebras; its gist being an identification of the
automorphism group of a Moufang loop with a subgroup of the automorphism group
of the associated group with triality.
\end{abstract}

%%%%%%%%%%%%%%%%%%%%%%%%%%%%%%%%%%%%%%%%%%%%%%%%%%%%%%%%%%%%%%%%%%%%%%%%%%%%%%%
% section INTRODUCTION                                                        %
%%%%%%%%%%%%%%%%%%%%%%%%%%%%%%%%%%%%%%%%%%%%%%%%%%%%%%%%%%%%%%%%%%%%%%%%%%%%%%%
\section{Introduction}

\noindent As we hope to attract the attention of both group- and
loop-theorists, we take the risk of being trivial at times and introduce most
of the background material carefully, although briefly. We refer the reader to
\cite{SprVel}, \cite{Pflugfelder}, \cite{Carter} and \cite{HallNagy} for a more
systematic exposition.

A groupoid $Q$ is a \emph{quasigroup} if the equation $xy=z$ has a unique
solution in $Q$ whenever two of the three elements $x$, $y$, $z\in Q$ are
known. A \emph{loop} is a quasigroup with a neutral element, denoted by
$\neutral$ in the sequel. \emph{Moufang loop} is a loop satisfying one of the
(equivalent) \emph{Moufang identities}, for instance the identity
$((xy)x)z=x(y(xz))$. The \emph{multiplication group} $\mlt{L}$ of a loop $L$ is
the group generated by all left and right translations $x\mapsto ax$, $x\mapsto
xa$, where $a\in L$.

Let $C$ be a vector space over a field $\field$, and $N:C\to F$ a
nondegenerate quadratic form. Define multiplication $\cdot$ on $C$ so that
$(C,\,+,\,\cdot)$ becomes a not necessarily associative ring. Then $C=(C,\,N)$
is a \emph{composition algebra} if $N(u\cdot v)=N(u)\cdot N(v)$ holds for
every $u$, $v\in C$. Composition algebras exist only in dimensions $1$, $2$,
$4$ and $8$, and we speak of an \emph{octonion algebra} when $\dim C=8$. A
composition algebra is called \emph{split} when it has nontrivial zero
divisors. By \cite[Theorem 1.8.1]{SprVel}, there is a unique split octonion
algebra $\octo{\field}$ over any field $\field$.

Write $\uocto{\field}$ for the set of all elements of unit norm in
$\octo{\field}$, and let $\sml{\field}$ be the quotient of $\uocto{\field}$ by
its center $\centre{\uocto{\field}}=\{\pm 1\}$. Since every composition
algebra satisfies all Moufang identities, both $\uocto{\field}$ and
$\sml{\field}$ are Moufang loops. Paige proved \cite{Paige} that
$\sml{\field}$ is nonassociative and simple (as a loop). Liebeck
\cite{Liebeck} used the classification of finite simple groups to conclude
that there are no other nonassociative finite simple Moufang loops besides
$\sml{\field}$, $\field$ finite.

Liebeck's proof relies heavily on results of Doro \cite{Doro}, that relate
Moufang loops to groups with triality. Before we define these groups, allow us
to say a few words about the (standard) notation. Let $G$ be a group. Working in
$G\rtimes \aut{G}$, when $g\in G$ and $\alpha\in\aut{G}$, we write $g^\alpha$
for the image of $g$ under $\alpha$, and $\comm{g}{\alpha}$ for
$g^{-1}g^{\alpha}$. Appealing to this convention, we say that $\alpha$
\emph{centralizes} $g$ if $g^\alpha=g$. Now, the pair $(G,\,S)$ is said to be a
\emph{group with triality} if $S\le \aut{G}$, $S=\spn{\sigma,\,\rho} \cong
S_3$, $\sigma$ is an involution, $\rho$ is of order $3$, $G=\comm{G}{S}$,
$\centre{GS}=\{1\}$, and the triality equation
\begin{displaymath}
    \comm{g}{\sigma}\comm{g}{\sigma}^\rho\comm{g}{\sigma}^{\rho^{2}}=1
\end{displaymath}
holds for every $g\in G$.

We now turn to geometrical loop theory. A \emph{$3$-net} is an incidence
structure $\net=(\points,\,\lines)$ with point set $\points$ and line set
$\lines$, where $\lines$ is a disjoint union of $3$ classes $\lines_i$ ($i=1$,
$2$, $3$) such that two distinct lines from the same class have no point in
common, and any two lines from distinct classes intersect in exactly one
point. A line from the class $\lines_i$ is usually referred to as an
\emph{$i$-line}. A permutation on $\points$ is a \emph{collineation} of $\net$
if it maps lines to lines. We speak of a \emph{direction preserving}
collineation if the line classes $\lines_i$ are invariant under the induced
permutation of lines.

There is a canonical correspondence between loops and $3$-nets. Any loop $L$
determines a $3$-net when we let $\points=L\times L$, $\lines_1=\{\{ (c,\,y)|y
\in L\}|c\in L\}$, $\lines_2=\{\{(x,\,c)|x\in L\}|c\in L\}$,
$\lines_3=\{\{(x,\,y)|x,\,y\in L$, $xy=c\}|c\in L\}$. Conversely, given a
$3$-net $\net=(\points,\,\lines)$ and the origin $\origin\in\points$, we can
introduce multiplication on the $1$-line $\ell$ through $\origin$ that turns
$\ell$ into a loop, called the \emph{coordinate loop} of $\net$. Since the
details of this construction are not essential for what follows, we omit them.

Let $\net$ be a $3$-net and $\ell_i\in\lines_i$, for some $i$. We define a
certain permutation $\bol{\ell_i}$ on the point set $\points$ (cf.\ Figure
\ref{Fg:Bol}). For $P\in\points$, let $a_j$ and $a_k$ be the lines through $P$
such that $a_j\in\lines_j$, $a_k\in\lines_k$, and
$\{i,\,j,\,k\}=\{1,\,2,\,3\}$. Then there are unique intersection points
$Q_j=a_j\cap \ell_i$, $Q_k=a_k\cap \ell_i$. We define $\bol{\ell_i}(P)=b_j\cap
b_k$, where $b_j$ is the unique $j$-line through $Q_k$, and $b_k$ the unique
$k$-line through $Q_j$. The permutation $\bol{\ell_i}$ is clearly an involution
satisfying $\bol{\ell_i}(\lines_j)=\lines_k$,
$\bol{\ell_i}(\lines_k)=\lines_j$. If it happens to be the case that
$\bol{\ell_i}$ is a collineation, we call it the \emph{Bol reflection with axis
$\ell_i$}.

%%%%% (FIGURE) [Fg:Bol] %%%%%%%%%%%%%%%%%%%%%%%%%%%%%%%%%%%%%%%%%%%%%%%%%%%%%%%
\setlength{\unitlength}{1mm}
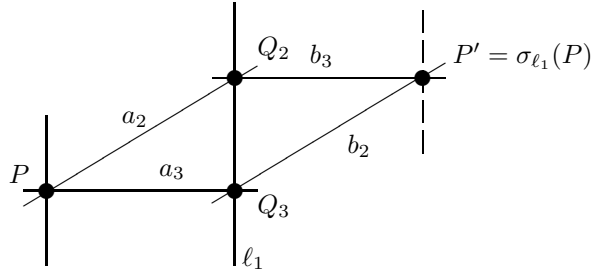
\begin{figure}
\begin{picture}(120,60)(-3,0)
\thicklines \put(50,10){\line(0,1){35}} \thinlines \put(25,10){\line(0,1){20}}
\multiput(75,25)(0,4){5}{\line(0,1){3}} \put(25,20){\circle*{2}}
\put(50,20){\circle*{2}} \put(50,35){\circle*{2}} \put(75,35){\circle*{2}}
\put(22,18){\line(5,3){31}} \put(78,37){\line(-5,-3){31}}
\put(22,20){\line(1,0){31}} \put(47,35){\line(1,0){31}} \put(20,21){$P$}
\put(79,37){$P^\prime = {\sigma_{\ell_1}}(P)$} \put(51,10){$\ell_1$}
\put(40,22){$a_3$} \put(35,29){$a_2$} \put(60,37){$b_3$} \put(65,25){$b_2$}
\put(53,17){$Q_3$} \put(53,38){$Q_2$}
\end{picture}
\caption{\label{Fg:Bol}The Bol reflection with axis $\ell_1$}
\end{figure}

It is clear that for any collineation $\gamma$ of $\net$ and any line $\ell$ we
have $\bol{\gamma(\ell)}=\gamma\bol{\ell}\gamma^{-1}$. Hence the set of Bol
reflections of $\net$ is invariant under conjugations by elements of the
collineation group $\coll{\net}$ of $\net$. A $3$-net $\net$ is called a
\emph{Moufang $3$-net} if $\bol{\ell}$ is a Bol reflection for every line
$\ell$. Bol proved that $\net$ is a Moufang $3$-net if and only if all
coordinate loops of $\net$ are Moufang (cf.\ \cite[p.\ 120]{Bruck}).

We are now coming to the crucial idea of this paper. For a Moufang $3$-net
$\net$ with origin $\origin$, denote by $\ell_i$ ($i=1$, $2$, $3$) the three
lines through $\origin$. As in \cite{HallNagy}, we write $\Gamma_0$ for the
subgroup of $\coll{N}$ generated by all Bol reflections of $\net$, and
$\Gamma$ for the direction preserving part of $\Gamma_0$. Also, let $S$ be the
subgroup generated by $\bol{\ell_1}$, $\bol{\ell_2}$ and $\bol{\ell_3}$.
According to \cite{HallNagy}, $\Gamma$ is a normal subgroup of index $6$ in
$\Gamma_0$, $\Gamma_0=\Gamma S$, and $(\Gamma,\,S)$ is a group with triality.
(Here, $S$ is understood as a subgroup of $\aut{\Gamma}$ by identifying
$\sigma\in S$ with the map $\tau\mapsto \sigma\tau\sigma^{-1}$.) We will
always fix $\sigma=\bol{\ell_1}$ and $\rho=\bol{\ell_1}\bol{\ell_2}$ in such a
situation, to obtain $S=\spn{\sigma,\,\rho}$ as in the definition of a group
with triality.

%%%%%%%%%%%%%%%%%%%%%%%%%%%%%%%%%%%%%%%%%%%%%%%%%%%%%%%%%%%%%%%%%%%%%%%%%%%%%%%
% section THE AUTOMORPHISMS                                                   %
%%%%%%%%%%%%%%%%%%%%%%%%%%%%%%%%%%%%%%%%%%%%%%%%%%%%%%%%%%%%%%%%%%%%%%%%%%%%%%%
\section{The Automorphisms}

\noindent Let $C$ be a composition algebra over $\field$. A map
$\alpha:C\to C$ is a \emph{linear automorphism} (resp.\
\emph{semilinear automorphism}) of $C$ if it is a bijective
$\field$-linear (resp.\ $\field$-semilinear) map preserving the
multiplication, i.e., satisfying $\alpha(uv)=\alpha(u)\alpha(v)$
for every $u$, $v\in C$. It is well known that the group of linear
automorphisms of $\octo{\field}$ is isomorphic to the Chevalley
group $G_2(\field)$, cf.\ \cite[Section 3]{Freudenthal},
\cite[Chapter 2]{SprVel}. The group of semilinear automorphisms of
$\octo{\field}$ is therefore isomorphic to
$G_2(\field)\rtimes\aut{\field}$.

Since every linear automorphism of a composition algebra is an isometry
\cite[Section 1.7]{SprVel}, it induces an automorphisms of the loop
$\sml{\field}$. By \cite[Theorem 3.3]{Vojtech}, every element of
$\octo{\field}$ is a sum of two elements of norm one. Consequently,
$\aut{\octo{\field}}\le\aut{\sml{\field}}$.

An automorphism $f\in\aut{\sml{\field}}$ will be called \emph{$($semi$)$linear}
if it is induced by a (semi)linear automorphism of $\octo{\field}$. By
considering extensions of automorphisms of $\sml{\field}$, it was proved in
\cite{Vojtech} that $\aut{\sml{\mathbb F_2}}$ is isomorphic to $G_2(\mathbb
F_2)$, where $\mathbb F_2$ is the two-element field. The aim of this paper is
to generalize this result (although using different techniques) and prove that
every automorphism of $\aut{\sml{\field}}$ is semilinear, provided $\field$ is
perfect. We reach this aim by identifying $\aut{\sml{\field}}$ with a certain
subgroup of the automorphism group of the group with triality associated with
$\sml{\field}$.

To begin with, we recall the geometrical characterization of automorphisms of
a loop.

%%%%% (LEMMA) [Lm:GeomChar] %%%%%%%%%%%%%%%%%%%%%%%%%%%%%%%%%%%%%%%%%%%%%%%%%%%
\begin{lemma}[Theorem 10.2 \cite{BarlStram}]\label{Lm:GeomChar}
Let $L$ be a loop and $\net$ its associated $3$-net. Any direction preserving
collineation which fixes the origin of $\net$ is of the form $(x,\,y)\mapsto
(x^\alpha,\,y^\alpha)$ for some $\alpha\in\aut{L}$. Conversely, the map
$\alpha:L\to L$ is an automorphism of $L$ if and only if $(x,\,y)\mapsto
(x^\alpha,\,y^\alpha)$ is a direction preserving collineation of $\net$.
\end{lemma}

We will denote the map $(x,\,y)\mapsto (x^\alpha,\,y^\alpha)$ by
$\sqmap{\alpha}$.

By \cite[Propositions 3.3 and 3.4]{HallNagy}, $\net$ is embedded in
$\Gamma_0=\Gamma S$ as follows. The lines of $\net$ correspond to the
conjugacy classes of $\sigma$ in $\Gamma_0$, two lines are parallel if and
only if the corresponding involutions are $\Gamma$-conjugate, and three
pairwise non-parallel lines have a point in common if and only if they
generate a subgroup isomorphic to $S_3$. In particular, the three lines
through the origin of $\net$ correspond to the three involutions of $S$.

As the set of Bol reflections of $\net$ is invariant under conjugations by
colli\-neations, every element $\varphi\in\coll{\net}$ normalizes the group
$\Gamma$ and induces an automorphism $\indmap{\varphi}$ of $\Gamma$. It is not
difficult to see that $\varphi$ fixes the three lines through the origin of
$\net$ if and only if $\indmap{\varphi}$ centralizes (the involutions of) $S$.

%%%%% (PROPOSITION) [Pr:Centralizer] %%%%%%%%%%%%%%%%%%%%%%%%%%%%%%%%%%%%%%%%%%
\begin{proposition}\label{Pr:Centralizer}
Let $L$ be a Moufang loop and $\net$ its associated $3$-net. Let $\Gamma_0$ be
the group of collineations generated by the Bol reflections of $\net$,
$\Gamma$ the direction preserving part of $\Gamma_0$, and $S\cong S_3$ the
group generated by the Bol reflections whose axis contains the origin of
$\net$. Then $\aut{L}\cong\centralizer{\aut{\Gamma}}{S}$.
\end{proposition}

\begin{proof}
Pick $\alpha\in\aut{L}$, and let $\indmap{\sqmap{\alpha}}$ be the automorphism
of $\Gamma$ induced by the collineation $\sqmap{\alpha}$. As $\sqmap{\alpha}$
fixes the three lines through the origin, $\indmap{\sqmap{\alpha}}$ belongs to
$\centralizer{\aut{\Gamma}}{S}$.

Conversely, an element $\psi\in\centralizer{\aut{\Gamma}}{S}$ normalizes the
conjugacy class of $\sigma$ in $\Gamma S$ and preserves the incidence
structure defined by the embedding of $\net$. This means that
$\psi=\indmap{\varphi}$ for some collineation $\varphi\in\coll{\net}$. Now,
$\psi$ centralizes $S$, therefore $\varphi$ fixes the three lines through the
origin. Thus $\varphi$ must be direction preserving, and there is
$\alpha\in\aut{L}$ such that $\varphi=\sqmap{\alpha}$, by Lemma
\ref{Lm:GeomChar}.
\end{proof}

It remains to add the last ingredient---groups of Lie type.

%%%%% (THEOREM) [Th:Main] %%%%%%%%%%%%%%%%%%%%%%%%%%%%%%%%%%%%%%%%%%%%%%%%%%%%%
\begin{theorem}\label{Th:Main}
Let $\field$ be a perfect field. Then the automorphism group of the
nonassociative simple Moufang loop $\sml{\field}$ constructed over $\field$
is isomorphic to the semidirect product $G_2(\field)\rtimes\aut{\field}$.
Every automorphism of $\sml{\field}$ is induced by a semilinear automorphism
of the split octonion algebra $\octo{\field}$.
\end{theorem}

\begin{proof}
We fix a perfect field $\field$, and assume that all simple Moufang loops and
Lie groups mentioned below are constructed over $\field$.

The group with triality associated with $M$ turns out to be its multiplicative
group $\mlt{M}\cong D_4$, and the graph automorphisms of $D_4$ are exactly the
triality automorphisms of $M$ (cf. \cite{Freudenthal}, \cite{Doro}). To be more
precise, Freudenthal proved this for the reals and Doro for finite fields,
however they based their arguments only on the root system and parabolic
subgroups, and that is why their result is valid over any field.

By \cite{Freudenthal}, $\centralizer{D_4}{\sigma} = B_3$, and by \cite[Lemmas
4.9, 4.10 and 4.3]{Liebeck}, $\centralizer{D_4}{\rho} = G_2$. As $G_2<B_3$, by
\cite[p.\ 28]{Gorenstein}, we have $\centralizer{D_4}{S_3} = G_2$.

Since $\field$ is perfect, $\aut{D_4}$ is isomorphic to $\Delta \rtimes
(\aut{\field} \times S_3)$, by a result of Steinberg (cf.\ \cite[Chapter
12]{Carter}). Here, $\Delta$ is the group of the inner and diagonal
automorphisms of $D_4$, and $S_3$ is the group of graph automorphisms of $D_4$.
When $\chr{\field}=2$ then no diagonal automorphisms exist, and
$\Delta=\inn{D_4}$. When $\chr{\field}\neq 2$ then $S_3$ acts faithfully on
$\Delta / \inn{D_4} \cong C_2 \times C_2$. Hence, in any case,
$\centralizer{\Delta}{S_3} = \centralizer{D_4}{S_3}$. Moreover, for the field
and graph automorphisms commute, we have $\centralizer{\aut{D_4}}{S_3} =
\centralizer{D_4}{S_3} \rtimes \aut{\field}$.

We have proved $\aut{M}\cong G_2\rtimes \aut{\field}$. The last statement
follows from the fact that the group of linear automorphisms of the split
octonion algebra is isomorphic to $G_2$.
\end{proof}

One of the open questions in loop theory is to decide which groups can be
obtained as multiplication groups of loops. Thinking along these lines we ask:
\emph{Which groups can be obtained as automorphism groups of loops}? Theorem
\ref{Th:Main} yields a partial answer. Namely, every Lie group of type $G_2$
over a perfect field can be obtained in this way.

Finally, the former author asked the latter one at the Loops '99 conference
whether $\aut{\sml{\field}}$ is simple when $F$ is finite. We now know that
this happens if and only if $\field$ is a finite prime field of odd
characteristic.

%%%%%%%%%%%%%%%%%%%%%%%%%%%%%%%%%%%%%%%%%%%%%%%%%%%%%%%%%%%%%%%%%%%%%%%%%%%%%%%
% BIBLIOGRAPHY                                                                %
%%%%%%%%%%%%%%%%%%%%%%%%%%%%%%%%%%%%%%%%%%%%%%%%%%%%%%%%%%%%%%%%%%%%%%%%%%%%%%%
\bibliographystyle{plain}

\end{document}